\documentclass[preprint,10pt]{elsarticle}
\usepackage{amssymb}
\usepackage{hyperref}
\usepackage{amsmath,amsthm}
\usepackage{amssymb,xcolor}
\usepackage{a4wide}
\hypersetup{linkcolor=blue, colorlinks=true
,citecolor = red}
\newtheorem{thm}{Theorem}[section]
\newtheorem{Def}{Definition}[section]
\newtheorem{Rem}{Remark}[section]

\newtheorem{corl}{Corollary}[section]
\newtheorem{example}{Example}[section]
\newtheorem{lem}{Lemma}[section]
\date{}
\catcode`\@=11
\def\theequation{\@arabic{\c@section}.\@arabic{\c@equation}}
\catcode`\@=12
\def\R{{I\!\!R}}

\newcommand{\Om} {\Omega}
\newcommand{\pa} {\partial}
\newcommand{\be} {\begin{equation}}
\newcommand{\ee} {\end{equation}}
\newcommand{\bea} {\begin{eqnarray}}
\newcommand{\eea} {\end{eqnarray}}
\newcommand{\Bea} {\begin{eqnarray*}}
\newcommand{\Eea} {\end{eqnarray*}}

\newcommand{\al} {\alpha}
\newcommand{\ba} {\beta}
\newcommand{\de} {\delta}
\newcommand{\ga} {\gamma}
\newcommand{\Ga} {\Gamma}

\newcommand{\De} {\Delta}

\newcommand{\rar}{\rightarrow}

\newcommand{\noi} {\noindent}

\newcommand{\D}{\Delta^2}

\newcommand{\e}{\epsilon}
\numberwithin{equation}{section}

\allowdisplaybreaks
\begin{document}

\begin{frontmatter}



\title{Isolated Singularities of Polyharmonic Operator in Even Dimension}


\author[rd]{Dhanya Rajendran}
\ead{dhanya.tr@gmail.com}
\author[as]{Abhishek Sarkar \corref{cor1}}
\ead{abhishek@math.tifrbng.res.in}

\address[rd]{Department of Mathematics\\
Indian Institute of Science,
Bangalore-560012, Karnataka, India.}

\address[as]{TIFR Centre For Applicable Mathematics\\
Post Bag No. 6503, Sharda Nagar, Bangalore-560065, Karnataka, India.}
\cortext[cor1]{Corresponding author}
\begin{abstract} We consider the equation $\De^2 u=g(x,u) \geq 0$ in the sense of distribution in $\Om'=\Om\setminus \{0\} $  where $u$ and $ -\De u\geq 0.$ Then it
is known that $u$ solves $\De^2 u=g(x,u)+\alpha \de_0-\beta \De \de_0,$ for some non-negative constants $\alpha$ and $ \beta.$ 
In this paper we study the existence of singular solutions to $\De^2 u= a(x) f(u)+\alpha \de_0-\beta \De \de_0$ in 
a domain $\Om\subset \R^4,$ $ a$ is a non-negative measurable function in some Lebesgue space. If $\De^2 u=a(x)f(u)$ in $\Om',$ then we
find the growth of the nonlinearity $f$ that determines $\alpha$ and $\beta$ to be $0.$ In case when $\alpha=\beta =0,$ we will establish
regularity results when $f(t)\leq C e^{\ga t},$ for some $C, \ga>0.$ 
This paper extends the work of Soranzo (1997) where the author finds the barrier function in higher dimensions $(N\geq 5)$ 
with a specific weight function $a(x)=|x|^\sigma.$ Later we discuss its analogous generalization for the polyharmonic operator.
\end{abstract}


\end{frontmatter}

\section{Introduction}
Isolated singularities of elliptic operators are studied extensively, see for eg. \cite{BL},\cite{Lions}, \cite{Soranzo}, \cite{T1} and \cite{T2}.
In this paper we wish to address the following problem and the questions related to it for 
the biharmonic(polyharmonic) operator in $\R^4 (\R^{2m})$:-\\[2 mm]
\textit{Question:} If a non negative measurable function $u$ is known to solve a PDE in the sense of distribution in a punctured domain, then what can one say about the 
differential equation satisfied by $u$ in the entire domain?\\[2mm]
In \cite{BL}, Brezis and Lions answered this question for the Laplace operator with the assumption that 
$$0\leq -\De u=f(u) \mbox{ in } \Om\setminus\{0\} \,\,,\,\, u\geq 0 \,\, ,\,\,\, \liminf_{t\rar \infty}\frac{f(t)}{t}>-\infty\,\,\, ,\,\,\, \Om\subset \R^N .$$ 
With the above hypotheses it was proved that  both $u$ and $f(u)$ belong to $L^1(\Om),$ and satisfy $-\De u= f(u)+\alpha \de_0,$ for some $\alpha\geq 0.$
For the dimension $N\geq 3,$ P.L.Lions\cite{Lions} found a sharp condition on $f$ that determines whether $\alpha$ is zero or not in the previous expression.
In \cite{DPS}, the authors further extended the result for dimension $N=2$ by finding the minimal growth rate of the function $f$ which guranteed $\alpha$ to be $0.$

Taliaferro, in his series of papers (see for e.g. \cite{T1}, \cite{T2}, \cite{Ghergu}) studied the isolated singularities of non-linear elliptic 
inequalities. In \cite{T2} the author studied the asymptotic behaviour of the positive solution of the differential 
inequality 
\begin{equation}\label{no1} 
0\leq -\De u\leq f(u) 
\end{equation}
 in a punctured domain under various assumptions on $f.$ If $N\geq 3$ and the function $f$ has a "super-critical" growth as in Lions\cite{Lions},
 (i.e. $\lim_{t\rar \infty}\frac{f(t)}{t^{\frac{N}{N-2}}}=\infty,$ ) then there exists arbitrarily 'large solutions' of (\ref{no1}). When
$N= 2,$ it was proved that there exists a punctured neighborhood of the origin
such that (\ref{no1}) admits arbitrarily large solutions near the origin, provided that 
$\log f(t)$ has a superlinear growth at infinity. Moreover author characterizes the singularity at the
origin of all solutions $u$ of (\ref{no1}) when $\log f(t)$ has a sublinear growth. Later Taliaferro, Ghergu and Moradifam in \cite{Ghergu} generalized these results to
polyharmonic inequalities. \\
The study of the polyharmonic equations of the type $(-\De)^{m} u=h(x,u)$ is associated to splitting the equation into a non-linear coupled system 
involving Laplace operator alone. 
Orsina and Ponce\cite{OP} proved the existence of solutions to 
$$(1)\,\,\left\{
\begin{array}{rlll}
-\De u &=& f(u,v)+\mu & \mbox{ in } \Om,\\
-\De v &=& g(u,v)+\eta & \mbox{ in }\Om,\\
u&=v&=0 &\mbox{ on } \pa \Om.
 \end{array}\right. 
$$ 
with the assumption that the continuous functions $f$ and $g$ are non increasing in first and second variables respectively with $f(0,t)=g(s,0)=0.$
But here the authors assumed that $\mu$ and $\eta$ are diffusive measures and Dirac distribution is not a diffusive measure. 
Considerable amount of existence/non-existence results have been proved for the problem $(1)$ when $f$ is a function of $v$ alone and $g$ depends only on $u$
and $\mu,\eta$ are 
 Radon measures. For eg. see \cite{VY} where the authors assumed $f(u,v)=v^p, $  $g(u,v)=u^q$ and with non-homogenous boundary condition. In \cite{GY} 
 authors dealt with sign changing functions $f$ and $g,$ with a polynomial type growth at infinity and the measure $\mu$ and $\eta$ were assumed to be multiples of $\de_0.$
 
 Our paper is closely related to the work of Soranzo \cite{Soranzo} where author considers the equation:
 $$\De^2 u= |x|^\sigma u^p \,\, \mbox{ with } u>0, \,\, -\De u>0  \mbox{ in } \Om\subset\R^N , N\geq 4 \mbox{ and }  \sigma \in (-4,0).$$
 A complete description of the singularity was provided when $1<p<\frac{N+\sigma}{N-4}$ for $N\geq 5,$ or $1<p<\infty$ when 
 $N=4.$ In this work we prove that the results of Soranzo can be improved for the dimension $N=4$ by replacing $u^p$ by more general exponential type function.

 %Our paper is closely related to the work of Soranzo \cite{Soranzo} and \cite{Soranzo1}.
 % If $f(t)$ is known to be bounded by $t^p,$ for some $p>1,$ then we can have take $a(x)\in L^k(\Om),$ \textcolor{red}{$k>r.$}
 % This work extends the results in \cite{Soranzo} and \cite{Soranzo1} for the dimension $N=2m.$

  \section{Preliminaries}
We assume that $\Om$ is a bounded open set in $\R^N , N\geq 4$ with smooth boundary and $0\in \Om.$ We denote $\Om'$ to be $\Om\setminus\{0\}.$ In this 
section we discuss some of the well known results for biharmonic operator.
\begin{thm} 
\label{thm1.1}
 (Brezis - Lions \cite{BL}) Let $u \in
L_{loc}^1(\Om')$ be such that $ \Delta u \in L_{loc}^1(\Om')$  in
the sense  of distributions in $\Om' ,$ $ u\geq 0$ in $\Om$ such that 

\begin{center}
 $-\Delta u+a u \geq  g \, \mbox{ a.e in } \Om,$
  \end{center}
where $a$ is a positive constant and {$g\in L_{loc}^1(\Om)$}. Then
there exist $\varphi \in L_{loc}^1(\Om)$ and $\alpha \geq$ 0 such
that
\begin{equation}\label{eqn1.1}
-\Delta u=\varphi+ \alpha\delta_0 \mbox{ in
}\mathcal{D}'(B_R)\end{equation} \noindent where $\delta_0$ is the
Dirac mass at origin.In particular, $u\in M^p_{loc}(B_R)$
\footnote { $M^p_{loc}(B_R)$ denotes the $ Marcinkeiwicz$ $
space$} where $p=N/N-2$ when $N\geq 3$ and $1\leq p<\infty$ is
arbitrary when $N=2$.
\end{thm}
\begin{thm} (Weyl Lemma, Simader\cite{Si})
 Suppose $ G\subset \R^N$ be  open
and let $u\in L^1_{loc}(G)$ satisfies
$$ \int_{G}u \D \varphi dx=0 \mbox{ for all } \varphi \in C^\infty_c(G),\,\, i.e.\,\, \D u= 0 \mbox{ in } \mathcal{D}'(G)  .$$
Then there exists $\tilde{u}\in C^\infty(G)$ with $\D \tilde{u}=0$
and $u=\tilde{u}$ a.e in $G.$
\end{thm}
\begin{thm}(Weak maximum principle:) Let $u\in W^{4,r}(\Om)$ be a
solution of
$$\left\{\begin{array}{cll} \D u = f(x)  \geq 0 \,\,\,& in \,\,
\Om\\
u\geq 0, -\Delta u \geq 0 & on \,\, \pa \Om\end{array}\right.$$
Then we have $u\geq 0$ and $-\Delta u\geq 0$ in $\Om.$
\end{thm}
\noi Proof of maximum principle easily follows by splitting the
equation into a (coupled) system of second order PDE's say:
$w=-\Delta u$ and $-\Delta w=f$ with the corresponding boundary conditions. \\
Using similar ideas we can infact prove a maximum principle with
weaker assumptions on the the smoothness of $u,$ which is stated below:
\begin{thm}\label{comp}
Let $u,\Delta u \in L^1(\Om)$ and $\D u \geq 0$ in the sense of
distributions. Also assume that $u, \Delta u$ are continuous near
$\pa \Om$ and $u>0, -\Delta u>0$ near $\pa \Om.$ Then $u(x)\geq 0$
in $\Om.$
\end{thm}
 
\begin{Def} Fundamental solution of $\D$ is defined as a locally integrable function $\Phi$
in $\R^N$ for which $\D \Phi=\delta_0$ and precisely expressed as
$$ \Phi(x)= a_N \left\{\begin{array}{lll}
|x|^{4-N} & \mbox{if} & N\geq 5\\
\log \frac{5}{|x|} & \mbox{ if }& N=4\\
|x| & \mbox{ if }& N=3\\
 |x|^2 \log \frac{5}{|x|} & \mbox{ if }& N=2\end{array}\right.
 $$ for some constant $a_N>0.$
\end{Def}
% We make the following observation regarding $\Delta |x|^\alpha$
% and $\D |x|^\alpha$ which will be useful in our calculations. For
% $x\neq 0$ and $-n<\alpha<0,$
% \begin{equation}\label{eqn1}
% \Delta |x|^\alpha= \alpha (\alpha -2+N)|x|^{\alpha-2} >0 \mbox{
% whenever } \alpha < 2-N.
% \end{equation}
% \begin{equation}\label{eqn2}
% \D |x|^\alpha = K_{N,\alpha} |x|^{\alpha-4} >0 \mbox{ whenever }
% \alpha<2-N  \mbox{ or } \alpha > 4-N.
% \end{equation}
\begin{thm}\label{main} Suppose $g:\Om'\times [0,\infty) \rar \R^+$ be a measurable function and let $u,\, \Delta u$ and $\D u \in
L^1_{loc}(\Om').$ Let $\D u=g(x,u)$ in $\mathcal{D}'(\Om')$
with $u\geq 0$ and $-\Delta u\geq 0$ a.e in $ \Om'.$ Then $u ,
g(x,u)\in L^1_{loc}(\Om)$ and there exist a non-negative constants
$\alpha$ , $\beta$ such that $\D u=g(x,u)+\alpha \delta_0-\beta \Delta \delta_0$ in
$\mathcal{D}'(\Om).$ 
% Infact $u\in M^{\frac{N}{N-4}}_{loc}(\Om)$ if$N\geq 5$ and $u\in L^q_{loc}(\Om)$ for all $q\in [1,\infty)$ if $N=4.$
\end{thm}
\noi Proof: Let us write $w=-\Delta u.$ Then $-\Delta w= g(x,u)\geq 0$ in
$\mathcal{D}'(\Om')$ and also given that $w, g(x,u)\in L^1_{loc}(\Om').$
Now as a direct application of Brezis-Lions Theorem \ref{thm1.1},
we obtain 
\begin{equation}\label{eqn11}
 -\Delta w= g(x,u)+\alpha \de_0 \mbox{ for some }\alpha\geq 0 \end{equation}
and $w, g(x,u)\in L^1_{loc}(\Om).$
Since $-\Delta u=w \geq 0$ in $\Om'$ again by Theorem \ref{thm1.1} $u\in L^1_{loc}(\Om)$ and 
$$-\Delta u= w+\beta \de_0 \mbox{ for some } \beta\geq 0.$$
Now substituting $w=-\Delta u -\beta \de_0$ in $(\ref{eqn11})$ we get 
\begin{equation}
 \De^2 u= g(x,u)+\alpha \de_0-\beta \De \de_0.
\end{equation}
Extending $g(x,u)$ to be zero outside $\Om$ we get 
$\D(u-f(u)*\Phi-\alpha \Phi-\beta \Gamma)=0$ in $\mathcal{D}'(\Om).$ By Weyl's
lemma for biharmonic operators, there exists a biharmonic function
$h \in C^{\infty}(\Om)$ and
$$ u= g(x,u)*\Phi+\alpha \Phi +\beta \Ga + h \,\, a.e \,\, \mbox{ in } \,\, \Om.$$
Note that $\Ga(x)$ belongs to Marcinkeiwicz space $ M^{\frac{N}{N-2}}(\Om)$ when $N\geq 2.$ By the property of the convolution of an $L^1$ 
function with the functions in $M^{\frac{N}{N-2}}(\R^N)$ we obtain $u\in
M_{loc}^{\frac{N}{N-2}}(\Om).$\hfill\qed.\\
The above result has been proved in \cite{Soranzo}(see Theorem 2) as an application of their main result on the system of equations. Proof is essentially
 based on the idea of Brezis-Lions type estimates. We have instead given a direct alternative proof for the same result.
% \begin{Rem}
%  This result has been proved in Soranzo(\cite{Soranzo}), Theorem 2 as an application of their main result on the system of equations which is essentially
%  based on the idea of Brezis-Lions estimates. We have instead given a direct alternative proof for the same result.
% \end{Rem}
\noi Theorem \ref{main} can be extended for polyharmonic operator in a standard way, for details see Theorem \ref{main1} .
% \begin{thm}\label{main2} Let $\Om\subset\R^{2m}$ and suppose $g:\Om'\times [0,\infty) \rar \R^+$ \textcolor{red}{is locally H$\ddot{o}$lder
% continuous (measurable is enough?)} and  $(\Delta)^i u \in
% L^1_{loc}(\Om')$ for $i=1,2\,\ldots m.$ Let $(-\De)^m u=g(x,u)$ in $\mathcal{D}'(\Om')$
% with $(-\De)^i u\geq 0$  a.e in $ \Om'$ for $1\leq i\leq m.$ Then $u ,\De u,\ldots \De^m u\in L^1_{loc}(\Om)$ and there exist a non-negative constants
% $\alpha_i$  such that $(-\De)^m u=g(x,u)+\alpha_1 \delta_0-\alpha_2 \Delta \delta_0+\ldots +(-1)^{m-1}\alpha_{m-1} \De^{m-1}\de_0$ in
% $\mathcal{D}'(\Om).$ 
% % Infact $u\in M^{\frac{N}{N-4}}_{loc}(\Om)$ if$N\geq 5$ and $u\in L^q_{loc}(\Om)$ for all $q\in [1,\infty)$ if $N=4.$
% 
%  \end{thm}

\section{Biharmonic operator in $\R^4$}
In this section we will restrict ourselves to the dimension $N=4$ and $g(x,u)$ to take a specific form $g(x,u)=a(x)f(u).$ Let $\Om $ be a bounded open set in
$\R^4,$ $ 0\in \Om$ and denote $\Om^{'}=\Om \setminus \{0\}$. 
We assume \\
\begin{itemize}
 \item [(H1)] $f:[0,\infty)\longrightarrow [0,\infty)$ is a 
continuous function which is nondecreasing in $\R^+$ and $f(0)=0.$
\item [(H2)] $a(x)$ is a non-negative measurable function in $L^k(\Om)$ for some $k>\frac{4}{3}.$
\item [(H3)] There exists $r_0>0$ such that ess$\inf_{B_{r_0}} a(x) >0.$\\
\end{itemize}

\par Let $u$ be a measurable function which solves the following problem: 
$$({P})\hspace{1cm}\left\{
\begin{array}{rllll}
\D u & = &a(x) f(u)\,\,\, \mbox{ in } \Om' \\
u\geq0 &, &-\Delta u\geq 0 \,\, \mbox{ in } \Om' 
% u \in C^4(\overline{\Om}\setminus\{0\}).
\end{array}\right.
$$
From Theorem \ref{main} we know that $u$ is a distributional solution of
$({P}_{\alpha,\beta})$
$$({P}_{\alpha,\beta})\hspace{1cm}\left\{
\begin{array}{lllll}\left.
\begin{array}{lllll}
\D u &= a(x) f(u)+\alpha \delta_0-\beta \De \de_0 \\
u\geq0& -\Delta u\geq 0 \end{array}\right\} \;\; \text{in} \;\Om, \\
\alpha,\beta  \geq 0,\, u \mbox{ and } a(x) f(u) \in L^1(\Om).
\end{array}\right.
$$\\[4 mm]
The assumption $(H3)$ suggests that the presence of such a weight function does not reduce the singularity of $a(x)f(u)$ at origin.
In particular, if $a(x)=|x|^\sigma$ for $\sigma\in (-3,0),$ then $a(x)$ satisfies $(H2)$ and $(H3).$\\[3 mm]
\noi Now assume that 
\begin{equation}\label{supqua}
 \displaystyle \lim_{t\rar \infty} \frac{f(t)}{t^2}=c\in (0,\infty].
\end{equation}
 i.e. 
$f(t)$ grows atleast at a rate of $t^2$ near infinity. Then for some $t_0$ large enough, we have 
$\displaystyle f(t)\geq \frac{c}{2}t^2$ for all $t\geq t_0.$ Suppose $u$ is a solution of $(P_{\al,\ba})$ and $f$
 satisfies \ref{supqua}. Then we know that for some biharmonic function $h$
$$u(x)= a(x) f(u)*\Phi+\alpha \Phi +\beta \Ga + h \,\, a.e \,\, \mbox{ in } \,\, \Om$$
where $\Phi$ is the fundamental solution of biharmonic operator in $\R^4$ and $\Ga$ is the fundamental solution of $-\De$ in $\R^4.$  Since 
$\al$ and $ a(x)f(u)$ are non-negative, we have $u(x)\geq \beta \, \Ga(x) + h(x).$
If $\beta\not = 0,$ fix  an $\tilde{r}\in (0,r_0)$ such that $\displaystyle u(x)\geq \frac{\beta}{2 \pi^2 |x|^2}\geq t_0$ whenever $|x|<\tilde{r}.$ Now, 
$$\int_{B_{\tilde{r}}}a(x) f(u) \geq C \int_{B_{\tilde{r}}}|x|^{-4} =\infty$$
which is a contradiction since $a(x)f(u)\in L^1(\Om).$ Thus $\beta=0$ if $f(t)$ grows at a rate faster than $t^2$ near infinity. We state this result in 
the next lemma.
\begin{lem}\label{lem1}
 Let $f$ satisfies the condition (\ref{supqua}) and $u$ solves $(P).$ Then for some $\alpha$ non-negative $\Delta^2 u = a(x) f(u)+\alpha \delta_0$ in 
 $\mathcal{D}'(\Om).$
\end{lem}

Now onwards we assume that $f$ satisfies (\ref{supqua}). We would like to address following questions in this paper: 
\begin{enumerate}
\item Can we find a sharp condition on $f$ that determines whether
$\alpha=0$ or not in $(P_{\alpha,0})$? \item If $\alpha=0$, is it true
that $u$ is regular in ${\Om} ?$
\end{enumerate}
\begin{Def}\label{def1}
We call $f$ a sub-exponential type function if $$\lim_{t \rar
\infty} f(t)e^{-\gamma t} \leq C \quad\rm{ for\quad some
\quad}\gamma, C>0.$$ We call  $f$ to be of super-exponential type
if it is not a sub-exponential type function.
\end{Def}
We will show that the above two questions can be answered based on
the non-linearity being a sub-exponential type function or not.
\begin{thm}(Removable Singularity)
Let $f$ be a super-exponential type function and $u$ is a distributional solution of 
$({P}).$ Then $u$ extends as a distributional solution of
$({P}_{0,0}).$
\end{thm}
\noi Proof: Given $u$ solves $(P),$ we know that $\D u= a(x)f(u)+\alpha \de_0-\beta \De \de_0$ for some
$\alpha,\beta \geq 0.$ To show the extendability of the distributional
solution we need to prove $\alpha=\beta= 0.$ Since $f$ is of super exponential type function, from Lemma \ref{lem1} it is clear that $\beta=0.$
 Let us assume that $\alpha>0$ and derive a contradiction. Note that we can find an $r$ small enough such that  
 $\displaystyle u(x)\geq -\frac{\alpha}{16 \pi^2}\log |x| $ whenever $|x|<r.$
 Since $f$ is not a sub-exponential type function, for a given $\ga>0$ there exists $t_0>0$ such that $f(t)\geq e^{\ga t}$ for all $ t\geq t_0.$ Thus, 
$$f(u(x))\,\geq\, f\left(-\frac{\alpha}{16 \pi^2}\log |x|\right)\, \geq\,\displaystyle e^{-\frac{\ga \alpha}{16 \pi^2}\log |x| }, \mbox{ for }|x|<<1 .$$
Now if we choose $\ga = \frac{64 \pi^2}{\alpha}$ in the above
inequality, it contradicts the fact $a(x) f(u) \in L^1(\Om).$ Thus
$\alpha=\beta = 0$ in $(P_{\al,\ba}).$ \hfill\qed\\
\begin{thm} {If $f(t)=t^p$ where $1 \leq p<\frac{4+\sigma}{2}$ and $a(x)=|x|^{\sigma},$ for $\sigma\in (-2,0)$ then  $(P_{\al,\ba})$ 
is solvable for $\al, \ba$  small enough.}
 \end{thm}
\noi Proof follows from Theorem 4(ii)) of Soranzo\cite{Soranzo}. The idea was to split the equation into a coupled system and find a sub and super solution
for the system. In the next theorem when $f$ satisfies \ref{supqua}, we find a super solution for $(P_{\al,0})$ directly without splitting the equation into
a coupled system and then use the idea of monotone
iteration to show the existence of a non-negative solution for $\alpha$ small enough. When $\beta\neq 0,$ such a direct monotone iteration technique
is not applicable as $\De \de_0$ is not a positive or a negative distribution, ie $\phi\geq 0,$ does not imply $\left<\De \de_0, \phi \right> \geq 0$ 
or $\leq 0.$
\begin{thm}
 Let $f$ and $a$ satisfy the hypotheses $(H1)-(H3).$ Additionally assume $\lim_{t\rar \infty}\frac{f(t)}{t^{2}}=c\in (0,\infty].$ Then 
 there exists an $\alpha_*>0$ such that for all $\alpha\leq \alpha_*$ the problem $(P_{\alpha,0})$ admits a solution in $B_r(0).$
\end{thm}
\noi Proof: We use the idea of sub and super solution to construct
a distributional solution for $(P_{\al,0})$ for $\al$ small enough.
Clearly $u_0=0$ is a subsolution for $(P_{\al,0}).$ Given that $f$ is a sub-exponential type nonlinearity,
 there exists a $\ga>0$ and a $C>0,$ such that $f(t)\leq C e^{\ga t}$ for all $t\in \R^+.$ 

Now define 
%\begin{equation}
 %\overline{u}(x)=\frac{-\log \ga |x|+C |x|^{3+\sigma}}{\ga} \mbox{ in } B_{r_1}(0).
%\end{equation}
\begin{equation}
 \overline{u}(x)=\frac{-\log |x|+C \phi}{\ga} \mbox{ in } B_{1}(0).
\end{equation}
where $\phi$ is the unique solution of the following Navier boundary value problem, 
\begin{eqnarray}\left\{ \begin{array}{rcl}
           \D \phi &=& \displaystyle -\frac{a(x)}{|x|}\log|x| \mbox{  in  } B_{1}(0)\\[3 mm]
           \displaystyle \phi=&0&=\De \phi \mbox{  on  } \pa B_{1}(0).             \label{supsoln}
          \end{array}\right.
          \end{eqnarray}
We notice that since $a(x)\in L^k(\Om),$ for some $k>\frac{4}{3},$ the term $ a(x)|x|^{-1}\log |x|\in L^p(B_1)$ for some $p>1.$ Hence the existence of a unique weak solution
$\phi\in W^{4,p}(B_1)$ is guarenteed by Gazzolla \cite{Gaz}, Theorem 2.20. Now by maximum principle we have $\phi\geq 0, -\De \phi\geq 0.$ \\
Therefore,% \begin{equation}
%  -\De \overline{u} = \frac{2}{\ga |x|^2}-\frac{C}{\ga}(3+\sigma)(5+\sigma)|x|^{1+\sigma}
% \end{equation}
\begin{equation}\label{ref11}
 \overline{u} \geq 0 \mbox{ in } B_1(0),
\end{equation}

\begin{equation}\label{ref22}
 -\De \overline{u} = \frac{2}{\ga |x|^2} - \frac{C}{\ga} \De \phi\, \geq\, 0.
\end{equation}
and 

\begin{equation}\label{ref33}
 \De^2 \overline{u}=\frac{\de_0}{8 \pi ^2 \ga}+\frac{C}{\ga |x|} {a(x)}\left|\log|x|\right|.
\end{equation}

\noi Note that $\displaystyle a(x) f(\overline{u}) \leq \frac{C}{|x|}a(x) e^{C \phi}.$  By Sobolev embedding, 
we know $W^{4,p}(\Om)\hookrightarrow C(\overline{\Om}),$ and hence $e^{C\phi}$ is bounded in $B_1(0).$ Now we fix an $r>0$ where  
$\displaystyle e^{C\phi}\leq \frac{|\log|x||}{\ga}$ in $B_r(0).$ We let $\Om=B_r(0)$ (where $r$ depends only on $\gamma$ and $C$) be a strict subdomian of $B_1(0)$ 
where $\displaystyle \frac{C}{\ga |x|} a(x)|\log|x|| \geq a(x) f(\overline{u}).$
Now from the choice of $r$ and equations \ref{ref11}) , (\ref{ref22}) 
and (\ref{ref33}) it is obvious that $\overline{u}$ is a super solution of $(P_{\al,0})$ where $\displaystyle \al=\frac{1}{8 \pi^2\ga}.$  Now let us define
inductively with $u_0=0$
$$
(P_{\al,0}^n)\,\,\,\left\{\begin{array}{rll}\displaystyle \D u_n=a(x)f(u_{n-1})+\al
\de_0 &
 \mbox{
in } \mathcal{D}'(\Om)\\
\displaystyle u_n>0 , -\Delta u_n>0 & \mbox{ in } \Om\\
\displaystyle u_n=\Delta u_n=0 &\mbox{ on }\pa \Om\\
%u_{n} \in M^{\frac{N}{N-4}}(\Om)
\end{array}   \right.
$$
Existence of such a sequence $\{u_n\}$ can be obtained by writing
$u_n=w_n +\al \Phi$ where
$$\left\{\begin{array}{lc}
\D w_n= a(x)f(u_{n-1}) \mbox{ in } \Om,\\
w_n=-\al \Phi , \Delta w_n= -\al \Delta \Phi \mbox{ on } \pa\Om, \\
w_n\in W^{4,r}(\Om) \mbox{ for some } r>1.\\
\end{array}\right.$$
Existence of $w_1$ is clear since $f(0)=0$ and from Theorem 2.2 of \cite{Gaz}. 
First let us show the positivity of $u_1$ and $-\Delta u_1$ in $\Om.$ Since $w_1$ is bounded, we can
choose $\e$ small enough so that $u_1=w_1+\al \Phi
>0$ and $-\De u_1>0$ in $B_\e.$ In $\Om \setminus B_\e$ by weak comparison
principle we can show that $u_1>0$ and $-\De u_1>0.$ Next we need to show that $u_1\leq \overline{u}.$ 
Note that by construction, $\overline{u}>0$ and $-\De \overline{u}>0$ in $\overline{B_r}\setminus\{0\}.$ Then,
$\overline{u}-u_1$ satisfies the set of equations 
$$\left\{  \begin{array}{rll}
   \D (\overline{u}-u_1)\geq 0 &\mbox{ in } \mathcal{D}(\Om),\\
    \overline{u}-u_1 >0 ,-\De(\overline{u}-u_1)>0&\mbox{ near }\pa\Om.
  \end{array}\right.
$$
Now using the maximum principle for distributional solutions (Theorem \ref{comp}) we find $u_1\leq \overline{u}.$

Assume that there exists a function $u_{k}$ solving $(P^k_{\al,0})$ for $k=1,2\cdots n$ and 
$$0\leq u_1\leq u_2 \ldots  \leq u_{n} \leq \overline{u} \mbox{ in }\Om.$$
Since $f$ is non-decreasing we have $a(x) f(u_{n})\in L^p(\Om),$ for some $p>1.$ Thus by Sobolev embedding there exists a 
$w_{n+1}\in C(\overline{\Om})\cap W^{4,p}(\Om).$ Also, 
$$\left\{\begin{array}{lc}
\D (u_{n+1}-u_n)=a(x)f(u_n)-a(x)f(u_{n-1}) \geq 0 &\mbox{ in } \Om\\
u_{n+1}=u_n , \Delta u_{n+1}=\Delta u_{n}& \mbox{ on } \pa\Om.
\end{array}\right.$$
Again from weak comparison principle $0<u_n\leq u_{n+1}$ and $0\leq -\De u_n\leq -\De u_{n+1}.$ As before one can show that $u_{n+1}\leq \overline{u}.$
Now if we define $u(x)= \lim_{n\rar
\infty}u_n(x)$ one can easily  verify that $u$ is a solution of $(P_{\al,0})$ for $\al=\frac{1}{8 \pi^2 \ga}.$
% \begin{equation}
%  \alpha_*=\sup \{{\alpha>0}: (P_{\alpha,0}) \mbox{ has a solution}\}
% \end{equation} 
For a given $f$ sub-exponential type function we have found a ball of radius $r$ such that $(P_{\alpha,0})$ posed on $B_r(0)$ has a
solution $u_\alpha$ for $\alpha=\frac{1}{8\pi^2\ga}.$ This solution $u_\alpha$ is a supersolution for $(P_{\alpha',0})$ posed in 
$B_r(0)$ and for $\alpha' \in (0,\alpha).$ Thus one can repeat the previous iteration and show that for all $\al' \in (0,\alpha)$
there exists a weak solution for $(P_{\al',0})$ in $B_r(0).$ % Thus $\alpha_*\geq \frac{1}{8\pi^2\ga}>0.$\textcolor{red}{ From the definition of $\alpha_*$ one can 
% infact check that for all $\al \in (0,)$ there exists a weak solution for $(P_{\al,0})$ in a ball around origin
% and the radius of the ball depends on $\alpha.$} please go through the proof and let me know if it is correct. does the radius of the ball shrink
% when $\alpha\rar 0?$ or we have a uniform ball for all $\alpha$ accordingly we will correct the corollory.
\hfill\qed
\begin{corl}
 Suppose for a given $\ga>0$ there exists a $C_\ga$ such that $f(t)\leq C_\ga e^{\ga t}$ for all $t\in \R^+.$ Then $(P_{\alpha,0})$ has a solution
 in $B_{r_\alpha}(0)$ for all $\alpha\in (0,\infty).$ In particular if $f(t)=t^p , p>2$ or $e^{t^\de}, \de<1$ then $(P_{\al,0})$ is solvable for all $\al>0.$
 \end{corl}
\noi Next we recall a Brezis-Merle \cite{BrezisMerle} type of
estimate for Biharmonic operator in $\R^4.$  Let $h$ be a
distributional solution of
$$(2)\,\,\,\,
\left \{\begin{array}{lll}
\D h= f & \mbox{ in } \Om\\
h=\Delta h=0 & \mbox{ on }\pa \Om.
\end{array}\right.
$$ where $\Om$ is a bounded domain in $\R^4.$ 
\begin{thm}(C.S Lin \cite{Lin})\label{Linthm}
Let $f\in L^1(\Om)$ and $h$ is a distributional solution of $(2).$ For a given $\de \in (0,32 \pi ^2)$ there exists
a constant $C_\de>0$ such that the following inequality holds:
$$\displaystyle \int_{\Om} \exp{( \frac{\de h}{\|f\|_{1}}}) dx \leq C_{\de} (\mbox{diam}\Om )^4$$
where diam $\Om$ denote the diameter of $\Om.$
\end{thm}
\begin{thm}
Let $f$ be a sub-exponential type function. Let $u$  be a solution 
of $({P}_{0,0})$ with $u=\De u= 0$ on $\pa \Om. $ Then $u$ is regular in $\Om.$
\end{thm}
\noi Proof: Let $u$ be a solution of $\D u= a(x) f(u) \mbox{ in } \Om$ with Navier boundary conditions. 
Write $g(x)=a(x) f(u),$ then $g \in L^1(\Om).$ Fix a $l>0$ and split $g=g_1+g_2$ where $\|g_1\|_1<\frac{1}{l}$ and $g_2\in L^\infty(\Om).$ Let 
$u_2$ be the unique solution of 
$$\left\{\begin{array}{lc}
\D u_2=g_2 \mbox{ in } \Om, \\
u_2=0, \Delta u_2=0 \mbox{ on } \pa\Om.
\end{array}\right.$$
Then 
$$\left\{\begin{array}{lc}
\D u_1=g_1 \mbox{ in } \Om, \\
u_1=0, \Delta u_1=0 \mbox{ on } \pa\Om.
\end{array}\right.$$
Choosing $\de=1$ in Theorem \ref{Linthm}, we find $\displaystyle \int_{\Om}exp(\frac{|u_1|}{\|g_1\|_1})< C_1 (diam \, \Om)^4. $ Thus $e^{l |u_1|}\in L^1(\Om).$
Since $u_2\in L^\infty(\Om),$ we have $e^{l |u|}\in L^1(\Om)$ for all $l>0.$ We use this higher intergrability property of $u$ in establishing its regularity. 

 We can  show that $a(x) f(u)\in L^r(\Om)$ for some $r>1.$
 In fact, 
$$\begin{array}{lll}
\displaystyle \int_{\Om}\left( a(x) f(u)\right)^r&\leq& \displaystyle \tilde{C}\int_{\Om} a(x)^r e^{\ga r u} \\
& \leq&\displaystyle  C_2 \left(\int_{\Om}a(x)^{ p r}\right)^{1/p}
\left(\int_{\Om}e^{p'\ga r u}\right)^{1/p'} <\infty\\
\end{array}$$
if we choose $p, r>1$ close enough to $1$ so that $1 < p.r\leq k,$ where $a(x)\in L^k(\Om).$ Now let 
$v$ be the unique weak solution of 
$$\left\{\begin{array}{lc}
\D v= a(x) f(u) \mbox{ in } \Om, \\
v=0, \Delta v=0 \mbox{ on } \pa\Om.
\end{array}\right.$$
 We have $v\in
C^{3,\gamma'}(\overline{\Om})$ for all $\gamma'\in (0,1).$ 
Now $u=v+h$ for some biharmonic function $h.$ Therefore $u\in C^{3,\ga'}(\Om).$
\hfill\qed.
\begin{Rem}
 The previous theorem is true even if $a(x)\in L^k(\Om)$ for some $k>1.$
\end{Rem}
\noi When $f$ is super exponential in nature an arbitrary solution of $\De^2 u= a(x)f(u)$ in $\mathcal{D}'(\Om)$ need not be bounded. We consider the following example.
\begin{example} Let $w(x)=(-4\log |x|)^{\frac{1}{\mu}}$ for some $\mu>1.$ Then one can verifty that whenever $x\not = 0,$
$$\De^2 w = b_1 e^{w^\mu} w^{1-4\mu}[b_2 w^{2 \mu}-b_3]$$
 for some positive constants $b_i.$ Since $f(w)= b_1 e^{w^\mu} w^{1-4\mu}[b_2 w^{2 \mu}-b_3]$ is super exponential in nature, $w$ extends as
 an unbounded distributional solution of $\De^2 w= f(w)$ in $B_r(0)$ for $r$ small enough.
\end{example}

\section{Polyharmonic Operator in $\mathbb{R}^{2m}$}
 We suppose $\Om$ is a bounded domain in $\mathbb{R}^N,\,\ N \geq 2m$ with smooth boundary and $0 \in \Om$. We denote $\Om'$ as $\Om \setminus \{0\}$.
\begin{thm}\label{main1} Suppose $g:\Om'\times [0,\infty) \rar \R^+$ is a measurable function and  $\De^{k}u \in
L^1_{loc}(\Om')$ for $k=0,1,..m$. Let $(-\De)^m u=g(x,u)$ in $\mathcal{D}'(\Om')$
with $(-\De)^{k}u\geq 0$ for $k=0,1,..,m-1$ a.e in $ \Om'.$ Then $u ,
g(x,u)\in L^1_{loc}(\Om)$ and there exist non-negative constants
$\alpha_0,...,\alpha_{m-1}$ such that $(-\De)^m u=g(x,u)+\displaystyle\sum_{i=0}^{m-1}{\alpha_i(-\De)^{i}} \delta_0$ in
$\mathcal{D}'(\Om).$ 
\end{thm}
%\textcolor{red}{write assumptions like H1, H2 ....}
Now we restrict ourselves to dimension $N=2m$ and $g(x,u)$ to take a specific form $g(x,u)=a(x)f(u)$. Throughout this section we make the following
assumption: \
\begin{itemize}
\item[$(H1')$] $f:[0,\infty) \mapsto [0,\infty)$ is a continuous function which is non-decreasing in $\mathbb{R}^+$ and $f(0)=0$.
\item[$(H2')$] $a(x)$ is non negative measurable function in $L^k(\Om)$ {for some $k>\frac{2m}{2m-1}$.}
\item[$(H3'$)] There exists $r_0>0$ such that $\mbox{essinf}_{B_{r_0}}a(x) >0$.\\
\end{itemize}
Let $u$ be a measurable function which satisfies the problem below,
$$({P^1})\hspace{1cm}\left\{
\begin{array}{cllll}
\begin{array}{rllll}
(-\De)^m u = a(x) f(u)\,\,\,& \mbox{ in } \Om' \\
(-\De)^k u\geq0 \mbox{ in } \Om', &k=0,..,m-1 \end{array}  \\
 u \in C^{2m}(\overline{\Om}\setminus\{0\}).
\end{array}\right.
$$
Then by \ref{main1} we know that $u$ is a distribution solution of $(P^1_{\alpha_0,..,\alpha_{m-1}})$ 
$$({P^1_{\alpha_0,..,\alpha_{m-1}}})\hspace{1cm}\left\{
\begin{array}{cllll}
\begin{array}{rllll}
(-\De)^m u &= a(x) f(u) +\displaystyle \sum_{i=0}^{m-1}\alpha_i(-\De)^{i}\delta_0 \text{ in } \Om \\
(-\De)^k u\geq0, &k=0,..,m-1 \text{ in } \Om'\end{array}  \\
 \alpha_i \geq 0, \text{ for }i=0,..,m-1 \text{ and } u, a(x)f(u) \in L^1(\Om).
\end{array}\right.
$$ 

In \cite{Soranzo1}, Soranzo et.al  considered a specific equation $(-\De)^m u= |x|^\sigma u^p$ in $\Om',$ with $\sigma \in (-2m,0)$ and $(-\De )^k u \geq 0,$ 
for $k=0,1, \ldots m.$ By Corollary 1 of \cite{Soranzo1}, if $N=2m$ and $p>\max\{1,\frac{N+\sigma}{2}\}$ then 
$\al_1 = \al_2 =\cdots =\al_{m-1}=0$ in $(P^1_{\al_0,\ldots \al_{m-1}}).$ This result can be sharpened for any weight function $a(x)$ satisfying 
$(H3)$ in a standard way and we skip the details of the proof.

% \begin{Rem}
% In particular if we take $a(x)= |x|^{\sigma}$ where, $-2m+1 <\sigma<0$. We note that $|x|^{\sigma} \in L^{\frac{2m}{2m-1}}$, and the above theorem holds.
% \end{Rem}

\begin{Rem}
Let $u$ satisfy $(P^1)$ and $\displaystyle \lim _{t \to \infty} \displaystyle \frac{f(t)}{t^m}=c \in (0,\infty].$ Then  
we have $ \alpha_1 = \alpha_2=..=\alpha_{m-1}=0$ in $(P^1_{\alpha_0,\ldots ,\al_{m-1}})$ and hence $u$ is a distributional solution of 
$(-\De)^m u= a(x)f(u)+\al_0 \de_0$ in $\Om.$
%te faster than $t^m$ near infinity. \textcolor{red}{do u wish to be more specific? you may quote the result
%exactly like in corollory 1 of \cite{Soranzo1} in our situation with weight function}
\end{Rem}
\par Now the following theorem gives us a sharp condition on $f$ which determines $\alpha_0 =0$ in $(P^1_{\al_0,0, \ldots , 0})$
and the proof is as similar to $Theorem \,\ 3.1$.
\begin{thm} Let $f$ be a super-exponential type function and $u$ is distribution solution of $(P^1)$.
Then $u$ extends as a distributional solution of $(P^1_{0,0,..,0})$.
\end{thm}

\begin{thm}
 Let $f$ and $a$ satisfy the hypotheses $(H1')-(H3').$ Additionally assume $\displaystyle \lim_{t\rar \infty}\frac{f(t)}{t^{m}}=c\in (0,\infty].$ Then 
 there exists an $\alpha_0>0$ such that for all $\alpha\leq \alpha_0$ the problem $(P^1_{\alpha,0,\ldots 0})$ admits a solution in $B_r(0),$ where
 the radius of the ball depends on the nonlinearity $f$.
\end{thm}
\noi{Proof:} We proceed as in Theorem 3.3, by constructing sub and super distributional solution for $(P^1_{\alpha,0,..,0})$ 
for all $\alpha$ small enough. We note that $u_0 =0$ is a sub-solution, and let 
\begin{equation}
 \overline{u}(x)=\frac{-\log |x|+C \phi}{\ga} \mbox{ in } B_{1}(0)
\end{equation}
where $\phi$ is the unique solution of the following Navier boundary value problem, 
\begin{eqnarray}\left\{ \begin{array}{rcl}
           (-\De)^m \phi &=& \displaystyle -\frac{a(x)}{|x|}\log|x| \mbox{  in  } B_{1}(0)\\[3 mm]
           \displaystyle \phi=\De \phi &=&..=(\De)^{m-1} \phi \mbox{  on  } \pa B_{1}(0).             \label{supsoln1}
          \end{array}\right.
          \end{eqnarray} 
          Then $\overline{u}$ is a supersolution of $(P^1_{\al,0\ldots 0})$ in a small ball $B_r(0)$. Rest of the proof follows exactly 
          as in the case of biharmonic operator.
 \hfill \qed

\par Next we state a Brezis-Merle type of type of estimates for poly-harmonic operator in $\mathbb{R}^{2m}$. 
\begin{thm} 
\label{thm1.1}
 (Martinazzi \cite{Mart}) Let $f \in L^1(B_R(x_0)), B_R(x_0) \subset \mathbb{R}^{2m},$ and let $v$ solve $$\left\{\begin{array}{cll} (-\De)^m v = f       \,\, \text{in} \,\, B_R(x_0),\\
 v= \D v=.....=\De^{m-1}v=0       \,\, \text{on } \pa B_R(x_0) \end{array}\right.$$ 
 Then, for any $p \in (0, \displaystyle\frac{\gamma_m}{\|f\|_{L^1(B_R(x_0))}}),$ we have $e^{2mp|v|} \in L^1(B_R(x_0))$ and $$\displaystyle \int_{B_R(x_0)}e^{2mp|v|} dx \leq C(p)R^{2m},$$
 where $\gamma_m= \displaystyle \frac{(2m-1)!}{2}\left|S^{2m}\right|$. 
  \end{thm}
  Finally with the help of above theorem we prove a regularity result for the polyharmonic operator.
  \begin{thm}
Let $a(x)$ and $f$ satisfies the properties as in $(H1')-(H3')$ and also assume that $f$ be a sub-exponential type function. Let $u$ be a solution $(P^1_{0,0,..,0})$ with $u= \De u =...=\De^{m-1}u=0$ on $\pa  \Om$.
Then $u \in C^{2m-1,\gamma'}(\Om)$, for all $\gamma' \in (0,1)$. %  is regular in $\Om$.\textcolor{red}{specify the regularity as we are not giving the proof, may be $C^{2m-1,\gamma'}$, check it.}
\end{thm}
\noi\textbf{Acknowledgement :}  Dhanya.R was supported by UGC under Dr.D.S Kothari Postdoctoral fellowship scheme No.F.4-2/2006(BSR)/13-1045. Both the authors
would like to thank Prof. S. Prashanth for various useful discussions.

\end{document}